\documentclass{article}
\usepackage{amsmath,amssymb,amsthm,amscd, mathtools, latexsym}
\usepackage{enumerate,varioref, dsfont, fancyhdr}
\usepackage{tikz-cd}
\usepackage{enumitem}
\usepackage{multicol}
\usepackage{hyperref}
\usepackage{url}
\usepackage{textcomp}

\newtheorem{Thm}{Theorem}[section]
\newtheorem{Lem}[Thm]{Lemma}
\newtheorem*{Lem*}{Lemma}

\newtheorem{Cor}[Thm]{Corollary}
\newtheorem{Conj}[Thm]{Conjecture}

\newtheorem{Def}[Thm]{Definition}
\newtheorem{claim}[Thm]{Claim}

{\theoremstyle{plain}

\newtheorem*{Ass*}{Assumption}

}

\def\pit{{\mathbb P}}

\def\0{{\mathcal O}}

\def\End{\mathop{\rm End}\nolimits}

\def\M{{\mathcal M}}

\usepackage[title]{appendix}

\begin{document}
\title{A Remark on the Rost Nilpotence Principle}
\author{Humberto Diaz}
\newcommand{\Addresses}{{\bigskip \footnotesize
\textsc{Department of Mathematics, Washington University, St. Louis, MO 63130} \par \nopagebreak
\textit{Email address}: \ \texttt{humberto@wustl.edu}}}

\date{}
\maketitle

\begin{abstract}
\noindent We consider a weak version of the Rost Nilpotence principle. For characteristic zero fields $k$, we show that if it holds for all smooth projective schemes over $k$, then the Rost Nilpotence Principle does also. We also correct the proof of a lemma in \cite{D}.
\end{abstract}
\section*{Introduction}
\noindent Let $X$ be a smooth projective scheme over a field $k$ of dimension $d$, by which we mean that all the irreducible components of $X$ are of dimension $d$. We will assume throughout that $k$ is perfect. As usual, for any field extension $F/k$, we use the notation $X_{F}:= X \times_{k} F$. Let $CH_{*} (-)$ (or $CH^{*}(-)$) denote the Chow group with integral coefficients indexed by dimension (resp., codimension) and let $\M_{k}$ denote the corresponding category of Chow motives. Finally, let $M(X) \in \M_{k}$ denote the Chow motive of $X$.\\
\indent For $\Gamma \in \End_{\M_{k}} (M(X))$, the Rost Nilpotence Principle predicts that if  there is some field extension $F/k$ for which the base extension $\Gamma_{F} = 0 \in \End_{\M_{F}} (M(X_{F}))$, then $\Gamma$ is nilpotent. This notion was introduced by Rost and proved in the case of quadrics over a field of characteristic $\neq 2$ \cite{R} (see also \cite{Br}). By work of \cite{CGM}, \cite{Gi1}, \cite{Gi}, \cite{Gi2}, \cite{RS} and \cite{D}, the Rost Nilpotence Principle also holds for surfaces, homogeneous varieties, as well as certain kinds of threefolds. It remains an open question whether the Rost Nilpotence Principle holds in general.\\
\indent In this short note, we consider the following weak version of this conjecture. 
\begin{Conj}[Weak Rost Nilpotence Principle]\label{WRNP} Let $\Gamma \in \End_{\M_{k}} (M(X))$. If  there is some field extension $F/k$ for which $\Gamma_{F} = 0 \in \End_{\M_{F}} (M(X_{F}))$, then 
\[\Gamma_{L*}: CH_{0} (X_{L}) \to CH_{0} (X_{L}) \]
is nilpotent for all field extensions $L/k$.
\end{Conj}  
\noindent Our main result below shows that under certain conditions, the Weak Rost Nilpotence principle actually implies the Rost Nilpotence Principle. Indeed, we have: 
\begin{Thm}\label{main} Suppose that $k$ is a field for which all schemes of dimension $\leq d-1$ admit a resolution of singularities and that the Rost Nilpotence Principle holds for all smooth projective schemes of dimension $\leq d-1$. If $\Gamma \in \End_{\M_{k}} (\M(X))$ satisfies Conjecture \ref{WRNP}, then $\Gamma$ satisfies the Rost Nilpotence Principle.
\end{Thm}
\noindent For the proof, we use a technique from a lemma in \cite{D}. As we note below, the proof of this lemma is not complete due to a gap. The author thanks Stefan Gille for pointing this out. We will use the main result to correct the proof. Besides this, we have the following curious consequence:
\begin{Cor}\label{last} Let $k$ have characteristic $0$. Then, Conjecture \ref{WRNP} holds for all smooth projective schemes of dimension $\leq d$ $\Leftrightarrow$ the Rost Nilpotence Principle holds for all smooth projective schemes of dimension $\leq d$.
\begin{proof} The ($\Leftarrow$) direction is immediate. The ($\Rightarrow$) direction can be proved by induction on $d$. The base step of $d=0$ is trivial. The induction step follows directly from Theorem \ref{main}. 
\end{proof}
\end{Cor}
\noindent It is perhaps surprising that (in characteristic zero, at least) proving that Rost Nilpotence holds in general can be reduced to proving Conjecture \ref{WRNP} holds in general, where the latter is only a statement about zero cycles. (See \cite{Gi2} for a similar type of result.) Applying Corollary \ref{last} to idempotents, one obtains the following:
\begin{Cor} Let $k$ have characteristic $0$ and $M \in \M_{k}$. The following are equivalent:
\begin{enumerate}[label=(\alph*)]
\item If $M_{F} = 0$ for some field extension $F/k$, then $CH_{0} (M_{L})=0$ for all field extensions $L/k$.
\item If $M_{F} = 0$ for some field extension $F/k$, then $M=0$.
\end{enumerate}
\end{Cor}
\section*{Proof of Theorem \ref{main}}
\noindent Let $\{X_{i}\}_{i=1}^{n}$ be the set of irreducible components of $X$ and $F_{i} = k(X_{i})$ the function field of $X_{i}$. So, assume that $\Gamma_{F} = 0 \in \End_{\M_{F}} (M(X_{F}))$ for some field extension $F/k$. Then, by the Weak Rost Nilpotence Principle, we have that 
\[\Gamma_{F_{i}*}: CH_{0} (X_{F_{i}}) \to CH_{0} (X_{F_{i}}) \]
is nilpotent. It follows that there is some $m > 0$ for which $\Gamma_{F_{i}*}^{m}$ acts trivially on $CH_{0} (X_{F_{i}})$ for each $i$. Thus, replacing $\Gamma$ with $\Gamma^{m}$, we can assume without loss of generality that $\Gamma_{F_{i}*}$ acts trivially on $CH_{0} (X_{F_{i}})$ for each $i$. Now, consider the following commutative diagram, where the (exact) rows come from the localization sequence:
\[\begin{tikzcd}
 & CH^{d} (X \times X) \arrow{d}{(\text{id}_{X}\times \Gamma)_{*}} \arrow{r} & \bigoplus_{i=1}^{n} CH_{0} (X_{F_{i}}) \arrow{d}{\oplus_{i} \Gamma_{F_{i}*}} \arrow{r} & 0\\
\bigoplus_{V} CH^{d-1} (V \times X)  \arrow{r} & CH^{d} (X \times X) \arrow{r} & \bigoplus_{i=1}^{n} CH_{0} (X_{F_{i}}) \arrow{r} & 0
\end{tikzcd}
\]
where the sum in the first term on the left is over the divisors of $X$. Since the rightmost vertical arrow is $0$, it follows that 
\[ \Gamma = (\text{id}_{X}\times \Gamma)_{*} (\Delta_{X}) \in CH^{d} (X \times X) \]
is supported on $V \times X$ for some divisor $V \subset X$. Then, we have the following lemma:
\begin{Lem}[\cite{D} Lemma A.2]\label{above} Suppose that resolution of singularities holds for schemes over $k$ of dimension $\leq r$. Let $V$ be a projective scheme over $k$ all of whose irreducible components have dimension $\leq r$. Then, there exists some smooth projective scheme $V'$ over $k$ all of whose irreducible components have dimension $\leq r$ and a morphism $\phi: V' \to V$ for which the push-forward
\[ (\text{id}_{X} \times \phi)_{*}: CH_{*} (X \times V') \to CH_{*} (X \times V) \]
is surjective. 
\end{Lem}
\noindent In fact, there exists $V'$ as in Lemma \ref{above} of dimension $=r$. To see this, note that for any irreducible component $V_{j}'$ of $V'$ of dimension $s \leq r$, one can replace $V_{j}'$ with $W_{j}= V_{j}' \times \pit^{r-s}$ and the morphism $V_{j}' \xrightarrow{\phi} V$ with the composition $W_{j} \xrightarrow{\pi_{V_{j}'}} V_{j}' \xrightarrow{\phi} V$; taking the union over all $j$, it is clear that the corresponding push-forward is then surjective.\\
\begin{claim}\label{claim} $\Gamma \in CH^{d} (X \times X) = \End_{\M_{k}} (M(X))$ is nilpotent.
\begin{proof}[Proof of Claim] Now, write $i: V \hookrightarrow X$, where $V$ is the subscheme of $X$ in the paragraph above Lemma \ref{above}. Then, let $V'$ be a smooth projective scheme of dimension $r=d-1$ as in the previous paragraph and denote the corresponding morphism $\phi: V' \to V$. As noted above, we have
\[ \Gamma = (\text{id}_{X} \times i)_{*}(\Gamma_{0}) \]
for some $\Gamma_{0} \in CH^{d-1} (X \times V)$. Let $\psi:= i \circ \phi: V' \to X$. Then, by Lemma \ref{above}, it follows that
\[ \Gamma = (\text{id}_{X} \times \psi)_{*}(\Gamma_{1}) = \Gamma_{\psi} \circ \Gamma_{1} \in CH^{d} (X \times X) \]
for some $\Gamma_{1} \in CH^{d-1} (X \times V')$. Following the argument of \cite{D} Lemma 2.8, we consider 
\[ \Phi =  \Gamma_{1} \circ \Gamma_{\psi} \in CH^{d-1} (V' \times V') = \End_{\M_{k}} (M(V'))  \]
Then, we have 
\[ \Phi^{2}_{F} = \Gamma_{1, F} \circ \Gamma_{F} \circ \Gamma_{\psi, F} = 0 \in CH^{d-1} (V'_{F} \times V'_{F}) = \End_{\M_{F}} (M(V'_{F})).\] 
Since the Rost Nilpotence Principle is assumed to hold for smooth projective schemes of dimension $\leq d-1$, it follows that $\Phi^{m} = 0$ for some $m$. Hence, 
\[ \Gamma^{m+1} = \Gamma_{\psi} \circ \Phi^{m} \circ \Gamma_{1} = 0 \in CH^{d} (X \times X) \]
as desired.  
\end{proof}
\end{claim}
\section*{A Correction}
\noindent In \cite{D}, the proof of Lemma 2.8 contains the incorrect statement that ``a sum of nilpotent correspondences is nilpotent" and, hence, the proof is not complete. This misstatement unfortunately also appears in Lemma A.1, but fortunately the only time Lemma A.1 is invoked in \cite{D} is in the proof of Lemma 2.8 (and, in fact, is not even needed then, as the proof below shows). As an application of Theorem \ref{main}, we offer the following completion of the proof. First, we recall the following terminology used in \cite{D}.
\begin{Def}\label{def} Given a scheme $X$ over a field $k$, we say that the Chow group of $X$ is {\em universally supported in dimension $\leq i$} if there exists a Zariski closed subset $V \subset X$, all of whose irreducible components are of dimension $\leq i$, for which the push-forward
\[ CH_{0} (V_{L}) \to CH_{0} (X_{L}) \]
is surjective for all field extensions $L/k$. 
\end{Def}
\begin{Lem}[\cite{D} Lemma 2.8]\label{supp} Assume that the following hold:
\begin{enumerate}[label=(\alph*)]
\item\label{Rost} the Rost nilpotence principle for irreducible smooth projective schemes of dimension $\leq d-1$;
\item resolution of singularities in dimension $\leq d-1$.
\end{enumerate}
Let $X$ be a smooth projective scheme of dimension $d$ whose Chow group is universally supported in dimension $\leq d-1$. Then, $X$ satisfies the Rost nilpotence principle.
\begin{proof}
By Theorem \ref{main}, it will suffice to prove that for any $\gamma \in \End_{\M_{k}} (M(X))$ for which $\gamma_{F}= 0$ for some extension $F/k$,  
\[ \gamma_{L*}: CH_{0} (X_{L}) \to CH_{0} (X_{L}) \]
is nilpotent for all field extensions $L/k$. As in Lemma 2.7 of \cite{D}, we may write
\[ \gamma = \gamma_{1} + \gamma_{2} \]
where $\gamma_{1}$ and $\gamma_{2}$ both vanish upon base extension to $F$ and where $\gamma_{1}$  is supported on $D \times X$ for $D \subset X$ a divisor and $\gamma_{2}$ is supported on $X \times V$ for $V \subset X$ a closed subscheme. By enlarging $V$ if necessary, we may assume that $V$ is also a divisor. By the discussion preceding Claim \ref{claim}, we have that there are smooth projective schemes of dimension $d-1$, $D'$ and $V'$, with morphisms $\phi_{D}:D' \to X$ and $\phi_{V}:V' \to X$ for which 
\begin{equation} \begin{split}
\gamma_{1} &= (\phi_{D} \times \text{id}_{X})_{*}(\gamma_{1}') = \gamma_{1}'\circ \prescript{t}{}{\Gamma_{\phi_{D}}},\\  
\gamma_{2} &= ( \text{id}_{X}\times \phi_{V})_{*}(\gamma_{2}') = \Gamma_{\phi_{V}}\circ \gamma_{2}'
\end{split}\label{layout}\end{equation}
for $\gamma_{1}' \in CH^{d-1} (D' \times X)$ and $\gamma_{2}' \in CH^{d-1} (X \times V')$. Here, the right hand sides of (\ref{layout}) follow from \cite{MNP} Lemma 2.1.3. From (\ref{layout}), it follows that
\[ \gamma_{1,L*}: CH_{0} (X_{L}) \to CH_{0} (X_{L}) \]
vanishes, since it factors through $\phi_{D}^{*}$. We deduce that
\begin{equation} \gamma_{L*}= \gamma_{2,L*} : CH_{0} (X_{L}) \to CH_{0} (X_{L}) \label{last-1}\end{equation}
Since $\gamma_{2}$ is supported on $X \times V$, $\prescript{t}{}{\gamma_{2}}$ is supported on $V \times X$. Then, Claim \ref{claim} applies, and it follows that $\prescript{t}{}{\gamma_{2}}$ is nilpotent and, hence, so is $\gamma_{2}$. So, (\ref{last-1}) is nilpotent, as desired.
\end{proof}
\end{Lem}

\end{document}